\documentclass[10pt]{article}
\usepackage{latexsym}
 \usepackage{amssymb}
\usepackage{amscd}

\newtheorem{theo}{THEOREM}[section]
\newtheorem{pro}[theo]{PROPOSITION}
\newtheorem{lemma}[theo] {LEMMA}

\newtheorem{rk}[theo]{Remark}
\newtheorem{defn}[theo]{DEFINITION}
\newtheorem{cor}[theo]{COROLLARY}
\newenvironment{proof}{\noindent \sf Proof. \sf}

\newcommand{\Mod}{\mathsf{Mod}}
\newcommand{\degree}{\mathsf{deg}}
\newcommand{\lo}{\longrightarrow}
\newcommand{\st}{\stackrel}
\newcommand{\be}{\beta}
\newcommand{\ga}{\gamma}
\newcommand{\cc}{{\mathcal{C}}}
\newcommand{\cb}{{\mathcal{B}}}

\newcommand{\de}{\delta}

\newcommand{\cm}{{\mathcal{M}}}

\newcommand{\cn}{{\mathcal{N}}}
\newcommand{\cGs}{{{\mathcal{C}}}_A \# G}

\newcommand{\ds}{\displaystyle}
\newcommand{\un}{\underline}

\newcommand{\ra}{\rightarrow}
\newcommand{\no}{\noindent}

\begin{document}

\sf

\title{Skew category, Galois covering and smash product of a $k$-category \footnote{18A32 16S35 16G20}}
\author{Claude Cibils,
\ \ Eduardo N. Marcos\thanks{ The second author wants to thank
CNPq (Brazil), for financial support, in form of a productivity
scholarship. We thank the IME of the Universidade de S\~{a}o Paulo for
support during the preparation of this work} }
\date{}
\maketitle
\begin{abstract}
In this paper we consider categories over a commutative ring
provided either with a free action or with a grading of a not
necessarily finite group. We define the smash product category and
the skew category and we show that these constructions agree with
the usual ones for algebras. In case of the smash product for an
infinite group our construction specialized for a ring agrees with
M. Beattie's construction of a ring with local units in \cite{be}.
We recover in a categorical generalized setting the Duality
Theorems of M. Cohen and S. Montgomery in \cite{cm}, and we
provide a unification with the results on coverings of quivers and
relations by E. Green in \cite{g}. We obtain a confirmation in a
quiver and relations free categorical setting that both
constructions are mutual inverses, namely the quotient of a free
action category and the smash product of a graded category.
Finally we describe functorial relations between the
representation theories of a category and of a Galois cover of it.
\end{abstract}

\section{Introduction}

Throughout this paper we will consider small categories $\cc$ over
a commutative ring $k$, which means that the objects $\cc_0$ forms
a set, the morphism set $_y\cc_x$ from an object $x$ to an object
$y$ is a $k$-module and the composition of morphisms is
$k$-bilinear.

We consider in addition $G$-categories and $G$-graded categories
over $k$ where $G$ is a group, inspired by the work of
\cite{cm,g,dlp,reri,mmm,cr} and by a preliminary version of
\cite{gms} . A $G$-category over $k$ is provided by a category
$\cc$ with firstly a set action of $G$ on the objects and secondly
$k$-module maps $s:\ _y\cc_x\rightarrow \ _{sy}\cc_{sx}$ for each
$s\in G$ and for each couple of objects $x$ and $y$, verifying
$s(gf)=(sg)(sf)$ in case $g$ and $f$ are morphisms which can be
composed in the category. Moreover for elements $t,s$ in the group
and a morphism $f$ we have $(ts)f= t(sf)$; we also require $1f=f$
where $1$ is the neutral element of $G$. In other words there is
group homomorphism from $G$ to the group of autofunctors of the
category. If the action of $G$ is free on $\cc_0$ we say that
$\cc$ is a free G-category.

A $G$-graded category $\cc$ over $k$ has a direct sum
decomposition of $k$-modules  $\ _y\cc_x\ = \oplus_{s\in G}\
{_y{\cc}_x}^s$ for each couple of objects, verifying ${\
_z\cc_y}^t {\ _y\cc_x}^s \subset {\ _y\cc_x}^{ts}$ and ${\ _x1_x}
\in   {\ _x\cc_x}^1 $.

In section \ref{galoiscoverings} we recall the definition of the
quotient category $\cc/G$ which makes sense only if $\cc$ is a
free $G$-category. It has been introduced and considered in
\cite{ri,boga} as well as in subsequent work by several authors
and recently in \cite{cr}. Note that by definition Galois
coverings are precisely given this way, namely $\cc$ is a Galois
covering of the quotient $\cc/ G$. Next we consider the
skew-category $\cc[G]$ which is an analogue of the skew-ring
construction. We prove a coherence result between our approach and
the ring theoretical construction in case $G$ is finite and the
category has a finite number of objects. The main result in this
section is that $\cc/G$ and $\cc[G]$ are equivalent.

In section \ref{smash} we consider a $G$-graded category $\cb$
over $k$ and we define the smash product category $\cb\#G$. In
case $G$ is finite, we obtain again a coherence result with the
$k$-algebra case used in \cite{cm}. The categorical approach for
the smash product that we introduce here is equivalent to the ring
without identity construction obtained by M. Beattie in \cite{be}.
Therefore the same remarks from the introduction of \cite{be} are
in force: the rings considered in \cite{mena,qu} differ from the
present approach. Note that in the context of quivers and
relations - that is categories given through a presentation - the
constructions described in \cite{g} corresponds to the smash
product introduced in this paper.

The first result in section \ref{smash} relates the categorical
smash product with quotients: let $\cb$ be a $G$-graded category
over $k$, then $\cb\#G$ is a free $G$-category and there is a
canonical isomorphism $\left(\cb\#G\right)/G\simeq\cb$. The second
result of this section concerns a free $G$-category $\cc$ over
$k$: the quotient category $\cc/G$ is graded and the smash product
$\left(\cc/G\right)\#G$ is isomorphic to $\cc$. As an immediate
Corollary we obtain that $(B\# G)[G]$ is equivalent to ${\cb}$. We
illustrate these results by considering the Kronecker category.

The combination of the quoted results provides a generalization of
the Cohen-Montgomery Duality Theorems \cite{cm} in the version
provided by M. Beattie \cite{be2}, for non necessarily finite
group and a category with possibly infinitely many objects.

In \cite{g} E. Green considers pairs $(\Gamma, \rho)$ formed by a
locally finite, path connected directed graph $\Gamma$ and a set
$\rho$ of relations, namely $k$-linear combinations of paths of
$\Gamma$. This data presents a category over $k$ : firstly
consider the free category on the set of paths - objects are
vertices and morphisms between two vertices are free $k$-modules
having a basis given by the paths between this vertices. Secondly
perform the quotient by the two-sided categorical ideal generated
by $\rho$. Regular coverings considered in \cite{g} concerns
categories over $k$ given by a fixed presentation which is locally
finite. We consider categories over $k$ with no need of a specific
presentation. In this context Galois coverings correspond to
regular ones, see for instance \cite{ri2}.

We provide in this paper a confirmation in an uniform pure
categorical setting that gradings of a category over $k$ are in
one-to-one correspondence with Galois coverings, after the work of
E. Green \cite{g}.

In order to make precise the link with \cite{cm} we  specialize
 our results to the case of a finite group $G$. Let $A$ be a
$G$-graded $k$-algebra, which we consider as a single object
category $\cc_A$ with endomorphism ring $A$. The smash product
algebra $A\# G$ has a free $G$-action and the corresponding skew
group algebra $(A\# G)[G]$ is proved in \cite{cm} to be a matrix
algebra over $A$ of size $|G|$. Now the smash product category
$\cc_A\# G$ we consider in this paper has $G$ as set of objects.
Since we assume here that $G$ is finite, we can consider the
$k$-algebra of the category $a(\cGs)$, namely the direct sum of
all the $k$-modules morphisms of $\cGs$ equipped with the usual
matrix product. We show that this algebra is precisely $A\#G$.

Our results in this specific case provide the existence of a
single object full subcategory of $\cGs[G]$ which is equivalent to
the ambient one and isomorphic to $\cc_A$. Since $\cGs$ has
precisely $|G|$ objects we infer Cohen-Montgomery result quoted
above.

In another direction note that Hochschild--Mitchell cohomology of
equivalent $k$-categories remain isomorphic through the
corresponding change of bimodules. In \cite{cr} a Cartan--Leray
spectral sequence is obtained relating the Hochschild--Mitchell
cohomology of the categories of a Galois covering
$\cc\rightarrow\cc/G$. Our results show that this Cartan--Leray
spectral sequence translates into a spectral sequence relating
Hochschild--Mitchel cohomology of $\cc$ and of $\cc[G]$, since
$\cc[G]$ is equivalent to $\cc/G$.

Finally in the last section we consider left modules over a
category $\cc$, namely contravariant functors from $\cc$ to the
category of $k$-modules. Note that if a category $\cc$ over $k$ is
presented by a quiver with relations, left $\cc$-modules are
precisely representations of the quiver subjected to the
relations. In case of a Galois covering of categories $\cc\ra\cb$,
modules over the base category $\cb$ coincide with the fixed
modules over the covering category under the action of the group
$G$, while all the modules over the covering category coincides
with graded ones over the base category. Note that in case of
rings, graded modules has been generalized to modules graded by a
$G$-set and general versions of the Duality Theorems of M. Cohen
and S. Montgomery has been obtained, see \cite{naraoy, nashoy,
rio}.

{\bf Acknowledgements:} The authors would like to thank Guillermo
Corti\~{n}as and Maria Julia Redondo for useful questions and remarks
preluding this work, as well as M. J. Redondo and A. Solotar for a
careful reading of a preliminary version of this paper. We also
thank the referee, in particular for useful comments concerning
previous results.

\section{Galois coverings and the
skew-category}\label{galoiscoverings}

Let $\cc$ be a free $G$-category over a commutative ring $k$.
Recall that the $G$-set of objects $\cc_0$ is free, namely if
$sx=x$ for an object $x$, then $s=1$.

\begin{defn} Given a free $G$-category $\cc$ over $k$, the quotient
category $\cc/G$ is the category over $k$ whose objects are the
$G$-orbits of $\cc_0$. For two orbits  $\alpha$ and $\beta$ the
$k$-module of morphisms from $\alpha$ to $\beta$ is
$$ _\beta(\cc/G)_\alpha = \left(  \bigoplus_{x\in\alpha,\  y\in \beta}{_y\cc_x}\right)/G.$$

\no Note that $\ds \oplus_{x \in \alpha,\ y \in \beta} \ _y \cc_x
$ is a $kG$-module and $\ds \left(\oplus_{x \in \alpha,\ y \in
\beta} \ _y\cc_x\right) /G$ is the largest $k$-module quotient
with trivial $G$-action: if $X$ is a $kG$-module $X/G=
X/(ker\epsilon)X$ where $\epsilon:kG\ra k$ is the augmentation
map, i.e. the $k$-linear map given by $\epsilon(s) = 1$ for all
$s\in G$.

Composition in $\cc/G$ is deduced from composition in $\cc$ and is
well defined precisely because the action of $G$ on the objects is
free.

By definition a Galois covering is the projection functor from a
free $G$-category to its quotient.
\end{defn}

\begin{lemma}\label{galois}
Let $\cc$ be a free $G$-category over $k$,  let $\alpha$ and
$\beta$ be orbits of objects and let $x_\alpha \in \alpha$ \ and
$x_\beta \in \beta$ be representatives. Then $\ds
_\beta(\cc/G)_\alpha$ and $\bigoplus_{s\in G}\ {
_{x{_\beta}}\cc_{sx_{\alpha}}}$ are canonically isomorphic.
\end{lemma}

\begin{proof}
Consider the normalization map $$\bigoplus_{x\in \alpha,\
y\in\beta} {_y}\cc{_x}\st{\phi}\lo \bigoplus _{s\in G}{
_{x_{\beta}}}\cc_{sx_{\alpha}}$$ given by $\phi(_yf_x) = sf $
where $s$ is the unique group element such that $sy = x_\beta.$
Note that $\phi(tf) = \phi(f)$ hence $\phi:\ _\beta(\cc/G)_\alpha
\lo \ds \oplus_{s\in G}\ {}_{x_{\beta}}\cc_{sx_{\alpha}}$ is well
defined. Conversely, consider
$$\psi:\bigoplus_{s\in G}\  _{x{_\beta}}\cc_{sx_{\alpha}}\hookrightarrow \bigoplus_{x\in\alpha, \ y\in\beta}{_y\cc_x\ra \ _\beta(\cc/G)_\alpha}$$ which is an inverse for $\phi$.

\end{proof}

\begin{defn}

Let $\cc$ be a $G$-category over $k$. The skew-category $\cc[G]$
has set of vertices
 $(\cc[G])_0 = \cc_0$, and morphisms $_y\left(\cc[G]\right)_x = \oplus_{s\in G}\  _y\cc_{sx}$.

Composition of morphisms is provided by the composition of $\cc$
after adjustment. In order to make this precise one needs to keep
track of the component where a morphism is located, we put $f=f_s$
in case $f\in {}_y\cc_{sx}$. This way $f_s\neq f_{s'}$ for $s\neq
s'$ in case $sx=s'x$. Let $g=g_t\in {_z}\cc_{ty}$ and $f=f_s\in
{_y}\cc_{sx}$. Then $gf = g \circ(tf)$.

\end{defn}

We show now a coherence result with the usual skew-algebra
construction.

\begin{pro}

Let $G$ be a finite group and let $\cc$ be a $G$-category over
$k$, with a finite number of objects. Let $a(\cc)$ be the
$k$-algebra associated to $\cc$, namely $\ds a(\cc) =
\oplus_{x,y\in\cc_0}\left({_y}\cc{_x}\right)$ provided  with the
matrix product induced by the composition of morphisms.

Then $G$ acts on $a(\cc)$ by algebra automorphisms and $a(\cc[G])$
is the skew-algebra $a(\cc)[G]$.

\end{pro}

\begin{proof}
Recall that if $A$ is a k-algebra with $G$ acting by automorphisms
on it, the skew group algebra $A[G]$ is the k-module $A\otimes kG$
with the twisted product provided by $sa =s(a)s$ for all $s\in G$
and $a\in A$. As it is well known and easy to check the twisting
(see for instance \cite{cascva}) conditions are satisfied and the
skew group algebra is a $k$-algebra which has $A$ and $kG$ as
subalgebras.

\no The isomorphism of algebras $\psi: a(\cc[G])\ra a(\cc)[G]$ is
defined as follows. Let $f_s\in\ _{y}\cc_{sx}$ be an elementary
morphism of $_y\cc[G]_x$. Then  $\psi(f_s) = f\otimes s$.

\no We have that $$\psi(gf) = \psi(g\circ(tf)) = g(tf)\otimes
ts,$$ while $$\psi(g)\psi(f) = (g\otimes t)(f\otimes s)=
g(tf)\otimes ts.$$

\no Moreover $$\psi(1) = \psi(\sum {_x}1_{x})= \sum \psi
({_x}1_{x})= \sum ({_x}1_{x}\otimes 1_G) = (\sum {_x}1_{x})\otimes
1_G = 1\otimes 1_G.$$
\end{proof}

The following result shows a fact that we need to be true in
$\cc[G]$. Together with the coherence property above it justifies
the definition we gave for the skew-category $\cc[G].$

\begin{lemma} Let $\cc$ be a $G$-category over $k$. Objects in the same $G$ orbit are isomorphic in $\cc[G].$

\end{lemma}

\begin{proof} If $sx = y$ then ${_y}1_y \in {_y}\cc[G]_x$ and ${_x}1_x\in {_x}\cc[G]_y$. These morphisms are mutual inverses.

 \end{proof}

\begin{pro} Let $\cc$ be a $G$-category over $k$ and let $\cc[G]$ be the skew-category. Choose an element $x_\alpha$ in each orbit $\alpha$ and let $\underline{\cc[G]}$ be the full subcategory of $\cc[G]$ given by these objects. Then the categories $\cc[G]$ and $\underline{\cc[G]}$ are equivalent.

\end{pro}

\begin{proof} Consider the functor $F:\cc[G]\ra \un{\cc[G]}$ defined on objects by
$F(x) = x_\alpha$ if the orbit of $x$ is $\alpha$. Since $x$ and
$x_\alpha$ are canonically isomorphic - see the Lemma above - a
morphism $f:x\ra y$ provides a unique morphism $x_\alpha \ra
x_\beta$ where $\beta$ is the orbit of $y$. This construction is
functorial, moreover $F$ and the natural embedding are inverse
equivalences.

 \end{proof}

\begin{pro}

Let $\cc$ be a free $G$-action category over $k$. Then the
categories $\cc/G$ and $\un{\cc[G]}$ are isomorphic.

\end{pro}

\begin{proof} This is a consequence of Lemma \ref{galois}.

 \end{proof}

We have proved the following

\begin{theo}
Let $\cc$ be a free $G$-category over $k$. The quotient category
$\cc/G$ and the skew category $\cc[G]$ are equivalent.
\end{theo}

\begin{rk}
This Theorem suggests that the categories $\cc[G]$ or
$\un{\cc[G]}$ can be considered as a substitute for the quotient
category in case the action of $G$ on $\cc$ is not free.
\end{rk}
\section{Smash Product Category}\label{smash}

Let $G$ be a group. A $G$-graded category $\cb$ over $k$ is a
category over $k$ together with a decomposition of each $k$-module
of morphisms ${_y}\cb_x = \bigoplus_{s\in G}\ {{_y}\cb_x}^s$ such
that $ {{_z}\cb_y}^t\ {_y\cb_x}^s \subset  \ {_z\cb_x}^{ts}$. In
particular ${_x}\cb_x$ is a $G$-graded algebra for each object
$x$.

The smash product category $\cb\# G$ is defined below. We will
show a coherence result in case $\cb$ is a single object category
given by a $G$-graded $k$-algebra. Then we will prove that a
Galois covering of categories is isomorphic to a smash product of
categories.

\begin{defn}
The smash product category $\cb\# G$ has object set $\cb_0\times
G.$ Let $(x,s)$ and $(y,t)$ be objects. The $k$-module of
morphisms is defined as follows: $$_{(y,t)}(\cb\# G)_{(x,s)} =
{{_y}\cb_x}^{t^{-1}s}.$$

Note that $_{(y,t)}(\cb\# G)_{(x,s)}$ and $_{(y,ut)}(\cb\#
G)_{(x,us)}$ are $k$-modules morphism from different objects which
coincide as $k$-modules.

In order to define the composition map
$$_{(z,u)}(\cb\# G)_{(y,t)}\ \otimes_k\  _{(y,t)}(\cb\#
G)_{(x,s)}\longrightarrow\ _{(z,u)}(\cb\# G)_{(x,s)}$$  note that
the left hand side is
 $$\ds _z\cb_y^{u^{-1}t}\ \otimes_{k}\ {{_y}\cb_x}^{t^{-1}x}$$ while the
right hand side is ${_z\cb_x}^{u^{-1}s}$. The graded composition
of $\cb$ provides the required map.

\end{defn}

\begin{pro}\label{quotientofsmash}
Let  $\cb$ be a $G$-graded category over $k$. The category $\cb\#
G$ is a free  $G$-category and $(\cb\# G)/G = \cb. $

\end{pro}

\begin{proof} On objects we define $u(x,s) = (x, us)$.

Let $f$ be a morphism in $ _{(y,t)}(\cb\# G)_{(x,s)} =\
_y\cb^{t^{-1}s}_x$, we define $$uf = f \in \ _{(y,ut)}\!(\cb\#
G)_{(x, us)}.$$

\no In other words the action is obtained by translation. Clearly
this is an action and moreover the action on objects is free. Note
that orbits of objects are in one to one correspondence with
objects of $\cb$, by retaining the first component of each couple.

This shows that $[(\cb\# G)/G]_0 \equiv \cb_0$. Moreover recall
that by definition we have

$$\ds _y[(\cb \# G)/G]_x = \left[\bigoplus_{s\in G, \ t\in G}\
_{(y,t)}(\cb\# G)_{(x,s)}\right]/G.$$

\no The action of $G$ on the numerator transports $_{(y,1)}(\cb\#
G)_{(x,s)}$ to $_{(y,u)}(\cb \# G)_{(x,us)}$ consequently those
$k$-modules are identified in the quotient. This remark shows the
following:
$$\ds \left[\bigoplus_{s\in G,\ t\in G}\ _{(y,t)}(\cb\# G)_{(x,s)}\right]/G =
 \bigoplus_{s\in G}\ _{(y,1)}(\cb\# G)_{(x,s)}
= \bigoplus_{s\in G}\ {{_y\cb}_x}^s\ .$$  \end{proof}

The following is a natural consequence of the former proposition.

\begin{cor} Let $G$ be a group and ${\cb}$ be a $G$-graded category over $k$.
Then the categories $(B\# G)[G]$ and ${\cb}$ are equivalent.

\end{cor}
\begin{proof}
We know that ${\cb}\#G$ is a Galois covering of ${\mathcal{B}}$.
The result of the previous section shows that $\cb[G]$ is
equivalent to $\cb/G$, namely $({\mathcal{B}}\#G)[G]$ is
equivalent to $({\mathcal{B}} \# G)/G = {\mathcal{B}}.$
\end{proof}

\begin{rk} In case a ring $A$ is graded by a finite group the coherence result
below shows that the preceding Corollary is Cohen-Montgomery
Duality Theorem, (see \cite{cm} and \cite{h}). Note that using
rings with local units, the Duality Theorem is still valid even if
$G$ is not finite by the results of \cite{be2,anma}, see also
\cite{ab}.
\end{rk}

\begin{rk}
In case of an infinite group grading of a $k$-algebra $A$ the
preceding Corollary shows that the category $(\cc_A\# G)[G]$ is
equivalent to $\cc_A$, where $\cc_A$ is the single object category
over $k$ determined by $A.$
\end{rk}

In order to prove the coherence property of the categorical smash
product with the smash product of a graded $G$-algebra we first
recall its definition. Note that the usual smash product of a
graded algebra requires $G$ to be finite, see for instance
\cite{cm}. As quoted in the Introduction, the categorical smash
product we consider in this paper is an alternative approach to
the ring without identity defined by M. Beattie in \cite{be}.

\begin{defn}

Let $A$ be a $G$-graded $k$-algebra where $G$ is a finite group.
The smash product $A\# k^G$ is the twisted (see for instance
\cite{cascva}) tensor product of algebras $A\otimes k^G$ where the
twist map $\tau:k^G\otimes A\ra A\otimes k^G$ is given as follows.
Let $\delta_s$ be the Dirac mass on an element $s\in G$, namely
$\delta_s$ the set map from $G$ to $k$ assigning $0$ to every
element except $s$, the value on $s$ is $1$. Let $f_t$ be an
element of degree $t$ in $A$. Then
$$\tau (\delta_t\otimes f_s) = f_s\otimes\delta_{s^{-1}t}.$$
\end{defn}

The verifications of the twisting properties of $\tau$ insuring
associativity are well known, we provide them for the convenience
of the reader and in order to avoid any missprint for $\tau$. We
omit the tensor product symbol as well as the map $\tau$, in other
words we consider $\de_tf_s = f_t\de_{s^{-1}t}$ as the product in
$A\otimes k^G.$

First observe that
$$(\de_u\de_t)f_s = 0 \mbox{ if
}  u\neq t \mbox{ while } (\de_u\de_t)f_s = f_s\de_{s^{-1}t}
 \mbox{ if } u = t$$ and $$\de_u(\de_t f_s) = \de_uf_s\de_{s^{-1}t}
= f_s\de_{s^{-1}u}\de_{s^{-1}t} =0 \mbox{ if } u\neq t  \mbox{
while } \de_u(\de_tf_s) = f_s\de_{s^{-1}t} \mbox{ if } u= t.$$

Second
 $$\de_u(g_tf_s) = g_tf_s \de_{(ts)^{-1}u}  \mbox{ since } g_tf_s\in A_{ts},$$
$$(\de_ug_t)f_s = g_t\de_{t^{-1}u}f_s = g_tf_s\de_{s^{-1}t^{-1}u} \ .$$
Moreover
$$\ds 1_{k^G}f_s = \sum_{u\in G}\de_uf_s = \sum_{s\in G}f_s\de_{s^{-1}u} = f_s1_{k^G}.$$

Recall that in case $\cb $ is a finite object set category over
$k$, the $k$-algebra of morphisms of $\cb$ is denoted  $a(\cb)$,
it is  the direct sum of the morphisms $k$-modules of $\cb$
equipped with the matrix product.

\begin{pro}
Let G be a finite group and let $A$ be a $G$-graded $k$-algebra.
Let $\cb_A$ be the single object category over $k$ with
endomorphism $k$-algebra $A$. Then the category $\cb_A$ is
$G$-graded and $a(\cb_A\# G) = A\# G$.
\end{pro}

Note that if $G$ is not finite $a(\cb_A\# G)$ can still be
considered (see \cite{be2, nashoy, rio}) as a ring with a set of
local units; the above proposition can also be proved in this more
general context where $A\# G$ is defined as a ring of this sort.

\begin{proof}
We define a map $\ds\phi:a(\cb_A\# G)\lo A\# G$ as follows.
Consider an elementary matrix $_tF_s\in a(\cb_A\#G)$ with all
entries zero except perhaps the entry $(x_0, t)-(x_0,s)$ with
value an element $f_{t^{-1}s}\in \ _{(x_0,t)}(\cb_A\#
G)_{(x_0,s)}$. Observe that
 $$_{(x_0,t)}(\cb_A\# G)_{(x_0,s)} = {\ _{x_0}(\cb_A)_{x_0}}^{t^{-1}s{}}= A_{t^{-1}s}.$$
We put $$\phi(_tF_s) = f_{t^{-1}s}\ \de_{s^{-1}}.$$

Note first that the unit element in $a(\cb\# G)$ is $\ds
1=\sum_{s\in G}{_s}1_s$ where $_s1_s$ is the elementary matrix
with nonzero entry only on the entry $(x_0, s)-(x_0,s)$ with value
$1\in A_{s^{-1}s}=A_1$. Then
$$\ds\phi(\sum_{s\in G}{_s}1_s) = \sum_{s\in G}1_1\ \de_{s^{-1}} = 1_{A\# G}.$$

Let $_vG_u$ and $_tF_s$ be elementary matrices. Then $_vG_u\ _tF_s
= 0$ if $u\neq t$ and  $_vG_u\ _uF_s $ is the elementary matrix
$_v(GF)_s$ with zero entries except a possible nonzero entry
$(x_0,v)-(x_0,s)$ with value $g_{v^{-1}t}f_{t^{-1}s}.$ We have
that
$$\phi(_vG_u)\phi(_tF_s)= g_{v^{-1}u}\ \de_{u^{-1}}f_{t^{-1}s}\ \de_{s^{-1}}
=g_{v^{-1}u}f_{t^{-1}s}\ \de_{s^{-1}tu^{-1}}\ \de_{s^{-1}}.$$ Note
that if  $u\neq t$ we have $\phi(_vG_u)\phi(_tF_s) = 0$ while if
$u=t$ $$\phi(_vG_u)\phi(_tF_s) =
g_{v^{-1}u}f_{t^{-1}s}\de_{s^{-1}} = \phi(_vG_t\ _tF_s). $$

Clearly the morphism we have defined is an injective and
surjective $k$-module map, since the image of the elementary
matrix $_tF_s$ belongs to $A_{t^{-1}s}\otimes\de_{s^{-1}}$.
\end{proof}

\medskip
The following result can be linked with Proposition
\ref{quotientofsmash}, in order to obtain that the categorical
smash product and the categorical quotient are inverse
constructions between graded categories and free action
categories.

\begin{theo}
Let $G$ be a group and let $\cc$ be a free $G$-category over $k$.
Then $\cc/G$ is $G$-graded and the category $\left(\cc/G\right)\#
G$ is isomorphic to $\cc$
\end{theo}
\begin{proof} Let $\{x_\alpha \mid x_\alpha\in \alpha\}$ be a choice of representatives in each $G$-orbit of objects of
$\cc$. Recall that $_\be(\cc/G)_\alpha =
\left[\ds\oplus_{x\in\alpha,\ y\in\be}\ {_y\cc_x}\right]/G$. Let
$_yf_x \in {_y\cc_x}$. In order to define the degree
$\mathsf{deg}[{_yf_x}]$ of the class $[{_yf_x}]$ we first
normalize $_yf_x$ in order to obtain as source object $x_\alpha$,
then we retain the gap between the target object obtained and
$x_\be$. More precisely there is an unique $s\in G$ such that
$sx=x_\alpha$, then there is a unique element of $G$ that we
denote $(\mathsf{deg}[f])^{-1}$ such that
$(\mathsf{deg}[f])^{-1}x_\be = sy$.

First we note that $\degree[f]$ is well defined: let $u\in G$,
then $uf\in{_{uy}\cc_{ux}}.$ The normalization uses $su^{-1}$ and
the target object is again $sy$.

Second we define $${_\be(\cc/G)_{\alpha}}^s = \left\{\sum[f]\
\mid\ f \in {_y\cc_x}\mbox{ where } x\in \alpha, \ y\in \be \mbox{
and } \degree[f] = s\right\}.$$ Then $_\be(\cc/G)_{\alpha}=  \ds
\bigoplus_{s\in
  G}{_\be(\cc/G)_{\alpha}}^s.$

Finally we check that $\degree(gf)=\degree(g)\degree(f)$. Indeed
consider $_{(\degree g)^{-1}x_\ga}g_{x_\be}$ and $_{(\degree
f)^{-1}x_\be}f_{x_\alpha}$. In order to compose those morphisms we
first have to adjust $g$ in order that it's  source object
coincide with the target object of $f$. The resulting composition
is
 $$_{(\degree f)^{-1}(\degree g)^{-1}x_\ga}\left[\left(\degree f^{-1}\right)gf\right]_{x_{\alpha}}.$$ Hence
$(\degree(gf))^{-1} = (\degree f)^{-1}(\degree g)^{-1}$ and
$\degree (gf) = (\degree g)(\degree f)$.

Using the graduation we have just described we define a functor
$$F:(\cc/G)\# G\ra \cc$$ as follows. Let $(\alpha, s)$ be an object
of the smash category. Then $F(\alpha, s) = sx_\alpha.$

Recall that by definition  $$_{(\be, t)}(\cc/G)_{(\alpha, s)} =
{_{\be}(\cc/G)_\alpha}^{t^{-1}s}.$$ Let $_{s^{-1}t
x_{\be}}[f]_{x_\alpha}\in \ _{(\be,
  t)}(\cc/G)\# G)_{(\alpha, s)}$  where
 $f\in\ _{s^{-1}tx_\beta}(\cc)_{x_\alpha}.$
We define $$F(_{s^{-1}t x_{\be}} [f]_{x_\alpha})\ =\  sf \ \in \
  _{tx_\be}\cc_{sx_{\alpha}}.$$

We consider now $$_{t^{-1}ux_\ga}[g]_{x_{\be}}\in \
_{(\ga,u)}(\cc\# G)_{(\be,t)}.$$ In order to compose $[g]$ and
$[f]$ in $(\cc/G)\# G$ we have to adjust $g$ using the $G$ action,
 obtaining the right object for its source namely the target of
$f.$

Then $[g] \circ  [f ]={ [s^{-1}tg]} [f]$. So
$$ F\left(_{s^{-1}ux_\alpha}[{s^{-1}tg}][f]_{x_\alpha}\right) = s\left((s^{-1}tg)
(f)\right)=(tg)(sf)= F(g)F(f)$$

It is not hard to see that $F$ is a bijection on objects and an
isomorphism of $k$-modules on morphisms.
\end{proof}

We end this section with a useful example. Let $\cb$ be the
Kronecker category with two objects $x$ and $y$, each object has
$k$ as endomorphism algebra, ${}_y\cb_x =ka\oplus kb$ while
${}_x\cb_y =0$. Let $G$ be a cyclic group with generator $t$. A
$G$-grading of $\cb$ is given by $\deg a = 1$ and $\deg b =t$. The
category $\cb \# G$ is the free category of the quiver having
vertices $(x,t^i)$ and $(y,t^j)$ and arrows from $(x,t^i)$ to
$(y,t^i)$ and to $(y,t^{i-1})$ for each $t^i\in G$. If $G$ is
finite this provides a crown, otherwise the quiver is infinite.
The quotient category is clearly the Kronecker category and the
quotient by the action of a proper subgroup provides an
intermediate crown.

In case $G$ is free abelian on two generators $s$ and $t$, the
grading given by $\deg a= s$ and $\deg b=t$ provides a non
connected Galois cover through the categorical smash product.

\section{Modules}
Let us recall the following definition:
\begin{defn} Let $\cc$ be a category over a commutative ring $k.$
A left $\cc$-module $\cm$ is a collection of $k$-modules
$\{_x\cm\}_{x\in \cc _0}$ provided with a left action of the
$k$-modules morphisms of $\cc$, given by $k$-module maps $\ds
_y\cc_x\otimes_k \ {_x\cm}\lo {_y\cm}$ where the image of
${_yf_x}\otimes {_xm}$ is denoted ${_yf_x}\ {_xm}$, verifying the
usual axioms
\begin{enumerate}
\item $_zf_y({_yf_x}\ {_xm}) =({_zf_y}\ {_yf_x}){_xm}$ \item
$_x1_x\ _xm = _xm$
\end{enumerate}
\end{defn}
In other words $\cm$ is a covariant functor from $\cc$ to the
category of $k$-modules. We denote by $\cc-\Mod$ the category of
left $\cc$-modules. As it is well known and easy to establish the
agreement property holds, namely if $\cc$ is a finite number
objects category and $a(\cc)$ is the corresponding algebra of
morphisms, $\cc$-modules and usual $a(\cc)$-modules coincide. In
particular if $A$ is a $k$-algebra and $\cc_A$ is the single
object category with endomorphism ring $A$, then $\cc_A$-modules
and $A$-modules are the same.

Let $F:\cc\ra\cb$ be a functor between categories over $k$. Note
that functors do not correspond really to algebra maps. In case
$\cc$ and $\cb$ have finite number of objects, the corresponding
map $a(\cc)\ra a(\cb )$ do not send, in general, $1_{a(\cc)}$ to
$1_{a(\cb)}$. However the map is additive and multiplicative.

\begin{defn}
Let $F:\cc\lo\cb$ be a functor of $k$-categories, and let $\cn$ be
a $\cb$-module. We denote by $F^*\cn$ the module $\cn\circ F$ and
by $F^*$ the corresponding functor $F^*:\cb-\Mod\lo \cc-\Mod$.
\end{defn}

In other words $F^*$ is the restriction functor and $F^*\cn$ is
the $\cc$-module given by $_x(F^*\cn) =\ _{F(x)}\cn$ and
$$(_yc_x).\left({_{F(x)}n}\right)\ = \ \left( _{F(y)}F(c)_{F(x)}\right)\ \left({_{F(x)}}n\right) \in\
_y(F^*\cn).$$

Let $f$ be a natural transformation between modules, then
$F^*(f)=f$. Note that $F^*$ is faithful. In case $F$ is a full
functor the functor $F^*$ is also full.

In case of a $G$-category $\cc$ over $k$ the group $G$ acts on the
category $\cc-\Mod$ in the following way: let $\cm$ be a
$\cc$-module and $s\in G$, then $^s\cm$ is given by $$_x(^{s}\cm)
= \ _{s^{-1}x}\cm$$ and $$\left({_yc_x}\right)
\left(_{s^{-1}x}m\right) = (s^{-1}c)m \in\ _{s^{-1}y}\cm = \
_y{(^s\cm)}.$$

\begin{theo}
Let $F:\cc\lo\cb$ be a Galois covering where $\cc$ is a free
$G$-category over $k$ and $\cb=\cc/G.$ Then $\cb-\Mod$ and
$(\cc-\Mod)^G$ are isomorphic categories.
\end{theo}

\begin{rk}
$(\cc-\mbox{Mod})^G$ is a full subcategory of $\cc$-Mod provided
by fixed objects under the action of $G$, namely modules $\cm$
such that $^s\cm =\cm$ for all $s\in G$.
\end{rk}
\begin{proof} The image of the functor $F^*$ is precisely $(\cc-\Mod)^G$. Indeed
$$_x\left(^s\left(F^*\cn\right)\right) =\ _{s^{-1}x}\left(F^*\cn\right) = _{F(s^{-1}x)}\cn = \
_{F(x)}\cn.$$

Note that by definition the functor $F$ of a Galois covering is
full.

\end{proof}

Next we want to prove that the entire category $\cc$-Mod is
isomorphic to the category of graded modules over $\cb$. Recall
that in case of a Galois covering $F:\cc\lo\cb$ (where $\cc$ is a
free $G$-category over $k$ and $\cb=\cc/G$) we have shown that
$\cb$ is $G$-graded. We also have shown that $\cc = \cb \# G$ and
we will use this description.

 First recall that in case $\cb$ is graded, a graded $\cb$-module is a $\cb$-module
 $\cn$ provided with a direct sum decomposition at each $k$-module, namely
$$_x\cn =\ds\bigoplus_{s\in G}{_x\cn^s}\mbox{ such that }_y\cb^t_x\ _x\cn^s\subset\  _y\cn^{ts}.$$
We denote by $\cb-\Mod_G$ the category of graded $\cb$-modules.
\begin{theo}
Let $\cb$ be a $G$-graded category over $k$ and let $\cc =\cb\# G$
be its Galois covering. Then $\cb-\Mod_G$ and $\cc-\Mod$ are
isomorphic categories.
\end{theo}

\begin{proof}
Let $\cn$ be a graded $\cb$-module. We define the $\cc$-module
$C(\cn)$ as follows:

$$_{(x,s)}{C(\cn)} = {}_x\cn^{s^{-1}},$$
and the action maps
$$_{(y,t)}\left(\cb\#
G\right)_{(x,s)}\otimes_k\ _{(x,s)}C(\cn) \lo \ _{(y,t)}{C(\cn)}$$
are given by the maps
$$_y\cb_x^{t^{-1}s}\otimes_k\  {_x\cn}^{s^{-1}} \lo\  _y\cn^{t^{-1}}.$$

Clearly this is an action of $\cc$ on $C(\cn)$ and a map between
graded $\cb$-modules gives in a functorial way a map between the
corresponding $\cc$-modules.

Conversely, let $\cm$ be a $\cc$-module. Let $B(\cm)$ be the
$\cb$-module defined by $$\ds _xB(\cm)= \bigoplus_{s\in
G}{_{(x,s)}}\cm.$$ The graduation is given by  $_xB(\cm)^s = \
_{(x, s^{-1})}\cm.$ The action is graded by construction. In order
to give the maps
$$_y\cb^t_x\otimes_k {_x}B(\cm)^s \lo\ _yB(\cm)^{ts}$$
recall that
$$_{(y,t)}(\cb\# G)_{(x,s)} = \ {_y\cb_x}^{t^{-1}s}.$$
Then
$${_y\cb_x}^t=\ _{\left(y,(ts)^{-1}\right)}{(\cb\#
G)_{(x,s^{-1})}}.$$ The required maps
$$_{\left(y,(ts)^{-1}\right)}(\cb\# G)_{(x, s^{-1})}\ {\otimes_k}\
_{(x,s^{-1})}{\cm}\lo
 \ _{\left(y,
(ts)^{-1}\right)}{\cm}$$ are obtained through the given action of
$\cb\# G$ on $\cm$.

Clearly the two constructions provide functors which are inverse
to each other.  \end{proof}

The following result is now clear:

\begin{pro}
Let $F:\cc\lo\cb$ be a Galois covering of categories over $k$ with
Galois group $G$. Let $\cn$ be a $\cb$-module and let $B$ be the
functor defined in the proof of the above Theorem. Then
$$BF^*(\cn)
=\ds \bigoplus_{s\in G}\cn^s \mbox{ where } \cn^s =\cn \mbox{ for
each } s\in G.$$ The graded action on $\ds\oplus_{s\in G}\cn^s$ is
given by ${_y\cb_x}^t\ {\otimes_k}\  {_x\cn}^s\ra {_y\cn}^t$ using
the $\cb$-action on $\cn$. Ignoring the grading, the $\cb$-module
$BF^*(N)$ is a direct sum of copies of $\cn$ indexed by $G$.
\end{pro}

\bibliographystyle{amsplain}

 \end{document}